\documentclass{article}
\usepackage{amsmath,amsthm,amsfonts}
\usepackage{amssymb}
\usepackage{epsfig}
\def\rank{\mathop{rank}}
\def\trace{\mathop{trace}}
\def\e{\epsilon}

\def\R{{\bf R}}

\def\id{\mathop{id}}
\newcommand{\N}{{\bf N}}
\usepackage[latin1]{inputenc}

\newtheorem{Def}{Definition}

\newtheorem{lemma}{Lemma}
\newtheorem{claim}{Claim}

\newtheorem{theorem}{Theorem}
\newtheorem{proposition}{Proposition}
\newtheorem{corollP}{Corollary}
\theoremstyle{remark}

\title{Almost solutions of equations in permutations}
\author{Lev Glebsky and Luis Manuel Rivera}

\begin{document}
\maketitle

\begin{abstract}
We will say that the permutations $f_1,...,f_n$ is an $\epsilon$-solution 
of an equation if the normalized Hamming distance between
its l.h.p. and r.h.p. is $\leq\epsilon$. We give a sufficient
conditions when near to an $\epsilon$-solution exists an exact
solution and some examples when there does not exist such a solution.



\end{abstract}
\section{Introduction and formulation of main results}

Let $S_n$ denote the group of all permutations of the finite set
$\{n\}=\{1,...,n\}$. For $ f,g\in S_n$ let  $h(f,g)=\frac{\mid
\{a:(a)f \neq (a)g\} \mid}{n}.$ It is easy to check that
$h(\cdot,\cdot)$ is a bi-invariant metric on $S_n$ \cite{MOP}. In
the article we are going to study almost solutions of equations in
$S_n$. For example, fix $p\in \N$ and for $f\in S_n$ consider
equation
\begin{equation} \label{eq1}
x^p=f.
\end{equation}
Such an $x$, if exists, is called to be $p$-root of $f$. Not every
$f\in S_n$ has a $p$-root, but any $f\in S_n$ has an almost $p$-root
for sufficiently large $n$. Precisely, the following theorem is true

\begin{theorem} \label{th1}
 For any $p\in\N$ there exists a sequence $\delta_n>0$, $lim_{n\to\infty}\delta_n = 0$
 such that for any $f\in S_n$ there exists a permutation $g \in S_n$
with  $h(g^p,f)\leq\delta_n$.
\end{theorem}

One of the motivations to consider such a question is studying sofic
groups, a class of groups that was introduced by Weiss and Gromov
\cite{w,grom}. And later G. Elek and  E. Szabó \cite{es1} define a
family of sofic groups that they called universal sofic groups with
the propertie that every sofic group is isomorphic to a subgroup of
an universal sofic group.

Theorem~\ref{th1} imply the following
\begin{corollP}
An universal sofic group $\cal U$ is an $\N$-root group. In other
words: any $g\in\cal U$ has a $p$-root for any $p\in\N$.
\end{corollP}
The other set of questions is about stability of a system of
equations in permutations, in order to give the precise formulation
of the problem we present the following definitions. Let
$w(x_1,...,x_k)$, $u(x_1,...,x_k)$,  be expressions using $x_j$,
$x_j^{-1}$ and multiplications (due to associativity we may think
that $w,u$ are words in $\{x_1,x_1^{-1},...,x_k,x_k^{-1}\}$).
\begin{Def}\label{def1}
\begin{enumerate}
\item We say that permutations $f_1,..,f_k$ are an $\e$-solution of an
equation ($\e$-satisfy an equation) $w(x_1,...,x_k)=u(x_1,...,x_k)$,
iff \\
$h(w(f_1,...,f_k),u(f_1,...,f_k))\leq\e$
\item We say that permutations $f_1,..,f_k$ are an $\e$-solution of a system
($\e$-satisfy a system) of equations
\begin{equation}\label{eq2}
w_i(x_1,...,x_k)=u_i(x_1,...,x_k),\;\; i=1,..,r
\end{equation}
 iff $f_1,..,f_k$
$\e$-satisfy every equation of the system.
\item  System (\ref{eq2}) is called stable (in permutations) iff there exists
$\delta_\e $, $\lim\limits_{\e\to 0}\delta_\e=0$ such that for any
$\e$-solution $f_1,f_2,...,f_k\in S_n$ of the system~(\ref{eq2})
there exists an exact solution $\tilde f_1,\tilde f_2,...,\tilde
f_k\in S_n$ of the system~(\ref{eq2}) such that $d(f_i,\tilde
f_i)<\delta_\e$ for $i=1,...,k$. (Note, that $\delta_\e$ is
independent of  $n$.)
\end{enumerate}
\end{Def}
There are some relations of stability of the system~(\ref{eq2}) and
the properties of the group $G=\langle
x_1,..,x_k\;|\;w_i(x_1,..,x_k)=u_i(x_1,...,x_k),\;i=1,...,r \rangle$.
\begin{theorem}\label{th2}
Let $G=\langle
x_1,..,x_k\;|\;w_i(x_1,..,x_k)=u_i(x_1,...,x_k),\;i=1,...,r \rangle$.\\
\begin{itemize}
\item If the group $G$ is finite then the system~(\ref{eq2}) is stable 
in permutations.
\item If the group $G$ is sofic but not residually finite, then the
system~(\ref{eq2}) is unstable in permutations.
\end{itemize}
\end{theorem}
So, for example the equation $x_1^3=x_2^{-1}x_1^2x_2$ is unstable in
permutations because the Baumslag-Soliter group $G=\langle
x_1,x_2\;|\;x_1^3=x_2^{-1}x_1^2x_2\rangle$ is sofic but not
residually finite \cite{Kropholler90, Irma, Radulescu02}. On the
other hand the system:
$$
x^3=y^3=(xy)^3=(x^2y)^3=\id
$$
is unstable in permutations, because the corresponding group is finite.
Of course, in most cases, Theorem~\ref{th2} says nothing about
stability of a system of equations and generally the question seems
to be very difficult. Particularly, we believe that the commutator
relation
\begin{equation}\label{eq3}
xy=yx
\end{equation}
is unstable but do not have a proof yet. In \cite{MR0148737} the similar but
easier question about commutator relation was considered.

The similar questions about matrices was widely studied and solved (
at least for commutator relation), see for example
\cite{MR928973,MR708811,MR1424963}. Let us discuss these results in
more details. First of all to formulate the problem we can
generalize Definition~\ref{def1} from $S_n$ to any family of sets,
where metrics and multiplications  are defined. Particularly, we may
ask $f_1,f_2,...,f_k$ in the definition~\ref{def1} to be unitary
(self-adjoint) matrices, with the  metrics $d(A,B)=\|A-B\|$,
$\|X\|=\sup_{\|x\|=1}\|Xx\|$. So, we can speak about stability of
the system~(\ref{eq2}) in
 unitary (self-adjoint) matrices. In this case it is also important
that $\delta_\e$ is independent of the size of the matrices.
The results of \cite{MR1424963,MR708811} say that the commutator relation is stable in
self-adjoint matrices and unstable in unitary matrices.

Although permutations have natural representations by unitary
matrices, instability of the  commutator relation for unitary
matrices seems to say nothing about stability of the commutator
relation in permutations.  One of the difficulty here is that the
representations of permutations by unitary matrices is not uniformly
continuous for the distances defined above. It looks like that the
following  distances for matrices are more relevant for the study
of stability in permutations.
\begin{enumerate}
\item $d(A,B)=\|A-B\|_T$, where $\|X\|_T=\sqrt{\frac{1}{n}\trace(XX^*)}$, or
\item $d(A,B)=\frac{\rank(A-B)}{n}$,
\end{enumerate}
where $n\times n $ is the size of the matrices. We do not know any
results about stability of commutator relations in matrices for
those distances, but there are some related works around von Neumann
algebras, where perturbations by compact operators is considered.
(Calkin algebras, essentially normal operators, see \cite{Farah07}
and the bibliography in it.)

\section{Proofs of the theorems}
In this section we present the proofs of theorems \ref{th1}  and \ref{th2}, in
order to proof theorem  \ref{th1} we use two propositions.\\
{\bf Proof of Theorem~\ref{th1}}

Some important facts. From the right and left invariance  of the
metric $d$ its follows that $d(x^n,y^n)\leq nd(x,y)$, the proof by
induction: $d(x^{(n+1)},y^{(n+1)})\leq
d(x^nx,x^ny)+d(x^ny,y^ny)=d(x,y)+d(x^n,y^n)$

So, it is enough to prove Theorem~\ref{th1} for prime $p$. Indeed,
if $d(f_1^{p_1},g)\leq \e_1$ and $d(f_2^{p_2},f_1)\leq \e_2$ then
$d(f_2^{p_1p_2},g)\leq \e_1+p_1\e_2$.

\begin{proposition}. \label{pr1} Let $f\in S_n$, let
$p$ be a prime number. The equation $x^p=f$ has an exact solution if
and only if for any $k\in\N$ the number of $kp$-cycles in $f$ is
divisible by $p$.
\end{proposition}
\begin{proof}

$\Rightarrow$  If  $x^p=f$, the $m$-cycles in $x$ with $(p,m)=1$
become $m$-cycles in
  $f$. For the $kp$-cycles in $x$, we obtain $p$ cycles of length $k$ in
  $f$. Therefore, the $kp$-cycles in $f$ can be obtained only by the
  $kp^2$-cycles in $x$, so, in $f$ we will have $p$ cycles of length 
$kp$ for every $kp^2$-cycles in $x$. 

  $\Leftarrow$ Given the permutation $f$, we will construct a permutation $x$
 that satisfies the equation. Suppose that the permutation $f$ has 
the following cyclic
 representation: $f=C_1 \dots C_h  D_1 \dots D_i$,
 where $C_i$ are all $i$-cycles, with $i$ relatively prime to $p$ and $D_k$ are
 all cycles of length $kp$. 
For any $C_r$ in $f$ we can write $C_r^{\alpha_r}$ in $x$ where
$\alpha_r p=1\mod r$.
Now, as for any $k$, the number $m$ of $kp$-cycles is divisible by  $p$, we can
divide the cycles $D_k$ in disjoint groups of size $p$, $D_k=d_1 d_2
...d_{m/p}$. For each group $d_l$
 $$
 d_l=(a_0^0,a_1^0,  \cdots , a_{kp-1}^0)(a_0^1,  a_1^1, \cdots,  a_{kp-1}^1)
\cdots (a_0^{p-1}, a_1^{p-1},  \cdots , a_{kp-1}^ {p-1})
$$
in $f$
 we can take
 $$x_l=(a_0^0,  a_0^1,  \cdots,  a_0^{p-1},  a_1^0,  a_1^1,  \cdots ,
a_1^{p-1},  a_2^0,  \cdots,  a_2^{p-1}, \cdots,  a_{kp-1}^0,
a_{kp-1}^1, \cdots , a_{kp-1}^{p-1}),$$ as the corresponding cycle
of the permutation $x$.
\end{proof}

\begin{proposition} \label{pr2}
Let $p$ be a prime number, let $f \in S_n$, then
there exist permutations $\tilde{f}, g \in S_n$, such that
$g^p=\tilde{f}$, and $h(\tilde{f},f)\leq
\frac{2\sqrt{2}(p-1)}{\sqrt{pn}}$.
\end{proposition}
\begin{proof}
In order to prove the proposition it is enough to construct
$\tilde{f}$ satisfying Proposition~\ref{pr1}.  Let the permutation
$f$ has the following cyclic
 re\-presentation: $f=C_1 \dots C_h  D_1 \dots D_j$, where the $C_i$ 
are all $i$-cycles, with $(p,i)=1$ and $D_i$ are all $ip$-cycles. Let
 $n_o$ be the number of all elements that belongs to the cycles $C_1,
 \cdots, C_h$, let $m_i$ be the number of all $ip$-cycles.  Because some of the
 $m_i$ can be zero, we consider the following set $S:=\{i \mid m_i \neq
 0\}$.\\
By Proposition~\ref{pr1}, in order to construct the permutation
$\tilde{f}$, we
 only need to change some cycles in $D_i$. We have $m_i=\alpha_i p +r_i$, $0\leq r_i<p$ and
 construct $\tilde{f}$ equal to $f$ but delete one element for the last $r_i$ $ip$-cycles, and make it
 fixed point ($(a_1,a_2,...,a_{ip-1},a_{ip})\to(a_1,a_2,...,a_{ip-1})(a_{ip})$).
 Then the distance between $f$ and  $\tilde{f}$ will be
$$
h(f,\tilde{f})= \frac{2 \sum_{i \in S}
r_i}{n}\leq\frac{2(p-1)|S|}{n},
$$
So, we only need to estimate $k=|S|$ for $n$ fixed. To make the
estimation let us put in order $S=\{s_1,s_2,...,s_k\}$, where $1\leq
s_1<s_2<...<s_k$. It follows that  $s_i\geq i$. Now
$$
n= n_0+p\sum_{i\in S}  m_i i\geq p\sum_{i=1}^k s_i\geq p
\sum_{i=1}^k i = p\frac{k(k+1)}{2}>p\frac{k^2}{2}.
$$
So, $\mid S \mid=k<\sqrt{2n/p}$ and the proposition follows.

\end{proof}

{\bf Proof of Theorem~\ref{th2} first part.}

Let  $V$ be a finite set of finite words in $x_1, x_2, ..., x_k$
that represent each element of the group $G$. Without loss of
generality we will assume that $\{x_1^{\pm},x_2^{\pm},...,x_n^{\pm}\}\subseteq V$. For
$v_1, v_2\in
  V$, the juxtaposed product $v_1 \cdot
v_2$ form a finite  word, which does not necessarily belong to $V$.
By
  the method of insertion and deletion of trivial and defining  relators of $G$,
  the word  $v_1 \cdot v_2$
  can be reduced to a word $v_{1,2} \in V$.
Let $m$ be the maximum length of the words appearing  during these
reduction processes for all triples of words $v_1, v_2, v_{1,2} \in
V$, $v_1v_2=v_{1,2}$ in $G$.

Let $f=\langle f_1,f_2,...,f_k\rangle\in S_n^k$ be an $\e$ solution
of System~\ref{eq2}. We think that the language of graphs is the
most appropriate to expose our proof. So, we can consider
$f_1,f_2,...,f_k\in S_n$  as an edge-colored graph $\Gamma(f)$ with
vertex set $V(\Gamma)=\{n\}$ and edge set $E(\Gamma)=E_1\cup
E_2\cup...\cup E_k$, where $E_i=\{(a,(a)f_i),a\in\{n\}\}$ is the
edges of color $i$. Let $N(a)$ be the $m$-neighborhood of a vertex
$a$ in $\Gamma$, where $m$ is the maximum length defined above. We
call $a\in\{n\}$ to be a good vertex iff for any $c\in N(a)$
 $f$ satisfies System~\ref{eq2} in $c$:
\begin{equation}\label{egv}
(c)w_i(f_1,...,f_k)=(c)u_i(f_1,...,f_k),\;\; i=1,..,r.
\end{equation}
A vertex is bad if it is not good.\\
\begin{claim}  Let $a\in \{n\}$ be a good vertex, then
$(c)x_i=(c)f_i$, for $c\in N(a)$ defines an action of $G=\langle
x_1,...,x_k\;|\;w_i(x_1,...,x_k)=u_i(x_1,...,x_k),\;i=1,...,r \rangle$
on $N(a)$. It implies that any $c\in N(a)$ is also a good vertex,
$N(c)=N(a)$ and
the set of good vertexes is disconnected from the set of bad vertexes.
\end{claim}\label{cl1}
\begin{proof}
Indeed, let $p_1=v_1v_2,p_2,...,p_n=v_{1,2}$ be the reduction from
$v_1v_2$ to $v_{1,2}$. Then, $(a)p_1(f)=(a)p_2(f)=...=(a)p_n(f)$ by
the definition of good vertexes and the claim follows.
\end{proof}

We may construct $(a)\tilde f_i=(a)f_i$ if $a$ is good vertex and
$(a)\tilde f_i=a$ if $a$ is bad vertex. It follows that $\tilde f_i$
satisfy System~\ref{eq2}, because the set of bad vertexes is
separated from the set of good vertexes.  So it is enough to show
that the set of bad vertexes is small. Let
$M=\{a\in\{n\}\;|\;au_i(f)\neq aw_i(f)\;\mbox{for some $i$}\}$, it is clear that
$|M|\leq \epsilon kn$. Then
the set of bad points is $M^*= \bigcup_{b \in M} N(b)$, so
$$
|M^*|\leq \sum_{b\in M} |N(b)|\leq \e kn
\left(1+k\frac{((2k-1)^{m-1}-1)}{k-1}\right)
$$
So $d(f_i,\tilde f_i)\leq C\e$, where $C$ depends only on the group
$G$.

{\bf Proof of the second part of Theorem~\ref{th2}}

For a group $X$ we will denote by $e_X$ the unity in $X$, some times we will write just 
$e$.
 For the second part of the theorem~\ref{th2} we recall the following
 definitions

\begin{Def}
Let $G$ be a group, $F \subseteq G$ be a finite subset, $\epsilon \geq
0$, and $\alpha >0$. An $(F, \epsilon, \alpha)$-representation in
$(S_n,h)$ is a map $\phi:F \rightarrow S_n$ with the following
properties

\begin{enumerate}
\item For any two elements $a,b \in F$, with $a \cdot b \in F$,
  $h(\phi(a)\phi(b),\phi(a\cdot b))< \epsilon$
\item If $e \in F$, then $\phi(e)=\id$
\item For any $a \neq e$, $h(\phi(a),\id)> \alpha$
\end{enumerate}
\end{Def}

\begin{Def}
The group $G$ is sofic if there exists $\alpha_0 >0$ such that for any
finite set $F \subseteq G$ and for any $\epsilon > 0$ there exists an
$(F, \epsilon, \alpha_0)$-representation in $(S_n,h)$.
\end{Def}

\begin{Def}
A group $G$ is residually finite iff for any $g\in G$, $g\neq e_G$, there exists a homomorphism 
$\phi$ to a finite group $H$ such that $\phi(g)\neq e_H$.
\end{Def}
We need the following lemma:
\begin{lemma}\label{l1}
If $d(\cdot,\cdot)$  is a bi-invariant metric, and $d(x_i,y_i)\leq \delta_i$,
i=1,...,r, then $d(x_1 \cdots x_r,y_1 \cdots y_r) \leq \sum_{i=1}^r \delta_i$
\end{lemma}
\begin{proof} By induction and by the bi-invariance of the metric, we have
  that

\begin{eqnarray*}\
d(x_1 \cdots x_{r+1},y_1 \cdots y_{r+1}) & \leq & d(x_1 \cdots
  x_{r+1},y_1 \cdots y_r \cdot x_{r+1}) \\[0mm]
&+&  d(y_1 \cdots y_r \cdot
  x_{r+1},y_1\cdots y_r \cdot y_{r+1} \\& \leq &   d(x_1 \cdots  x_{r},y_1
  \cdots y_r ) + d(x_{r+1},y_{r+1})  \\& \leq & \sum_{i=1}^r \delta_i +  \delta_{r+1}=\sum_{i=1}^{r+1}\delta_i\\
\end{eqnarray*}
\end{proof}
In order to proof the second part of theorem \ref{th2}, we will prove the following
\begin{proposition}
If  $G=\langle
x_1,\ldots,x_k\;|\;w_i(x_1,\ldots ,x_k)=u_i(x_1,\ldots ,x_k),\;i=1,\ldots,r \rangle$
is sofic and System~\ref{eq2} is stable, then $G$ is residually
finite.
\end{proposition}
\begin{proof}

Let $p(\bar{x})$ be any word in $G$, $p\neq e_G$. We need to
construct a homomorphism $\phi$ to a finite group, such that
$\phi(p)\neq e$. We denote by $V^*$ the set of the words
$w_i(x_1,\ldots ,x_k),\;u_i(x_1,\ldots,x_k),\;i=1,...,r$ and all its
subwords and let $p^*$ denote the set of all subwords of $p$. Let
$F$ the following set: $F:=\{1,x_1,x_1^{-1},\ldots ,x_k, x_k^{-1} \}
\cup V^* \cup p^*$. It is clear, that $F$ is finite and for
any word in $F$ all its subwords belongs to $F$. As the group $G$ is
sofic, there exists $\alpha
> 0$ such that for any $\epsilon > 0$ there exists an $(F,\epsilon,
\alpha)$-representation $\varphi$. We denote  $a_i:=\varphi(x_i)$.
\begin{claim}
For any word $v(\bar{x}) \in F$, $h\Big(\varphi\big(v(\bar{x})\big),v(\bar{a})\Big)
<(2 \mid v(\bar{x}) \mid -1)\epsilon$
\end{claim}\label{cl2}
\begin{proof}  As $h\big(\varphi(x_i^{-1}),a_i\big)< \epsilon$, by induction
and Lemma~\ref{l1} we have
\begin{eqnarray*}
h\Big(\varphi\big(x_i^{\pm 1}v(\bar{x})\big),a_i^{\pm 1}v(\bar{a})\Big)
  &<& h\Big(\varphi\big(x_i^{\pm 1}v(\bar{x})\big),\varphi\big(x_i^{\pm 1}\big)
\varphi\big(v(\bar{x})\big)\Big) \\[0mm]
&+& h\Big(\varphi\big(x_i^{\pm 1}\big)\varphi\big(v(\bar{a})\big),a_i^{\pm 1}
  v(\bar{a})\Big)\\ &\leq&
\epsilon + (\epsilon + (2\mid v(\bar{x}) \mid-1) \epsilon)=(2(1+\mid
  v(\bar{x})\mid )-1)\epsilon.
\end{eqnarray*}

\end{proof}
 So, $a_i=\varphi(x_i)$ is an
$\epsilon^*$-solution of the system \ref{eq2} with
$\epsilon^*=max\{2(\mid w_i \mid +
  \mid u_i \mid ) \epsilon\}$. As the  system \ref{eq2} is stable,
   we can find an exact solution $b_1, ... , b_n$ of the system
  \ref{eq2}, with $h(a_i,b_i) \leq \delta_{\epsilon ^*}$ for any $i$.
Then $\phi(x_i)=b_i$, can be extended to a homomorphism $G\to S_n$.

\begin{claim} For any word $v\in F$, $v\neq e_G$ one has
$h\Big(\phi\big(w(\bar{x})\big),id\Big) \geq \alpha - \delta^*_{\epsilon ^*}\mid
w(\bar{x}) \mid -2 \mid w(\bar{x})
  \mid \epsilon$
\end{claim}\label{cl3}

\begin{proof}
Because $\phi(x_i)$ is a homomorphism $\phi\big(w(\bar{x})\big)=w(\bar{b})$,
and by Lemma \ref{l1}, $h\big(w(\bar{a}),w(\bar{b})\big) \leq
\delta^*_{\epsilon ^*}|w(\bar{x})|$ then
\begin{eqnarray*}
h\Big(\phi\big(w(\bar{x})\big),\varphi\big(w(\bar{x})\big)\Big) &\leq&  
h\Big(\phi\big(w(\bar{x})\big),w(\bar{b})\Big)
 + h\Big(w(\bar{b}), \varphi\big(w(\bar{x})\big)\Big)\\ &\leq&
h\Big(w(\bar{b}), w(\bar{a})\Big) + h\Big(w(\bar{a}),\varphi\big(w(\bar{x})\big)\Big)\\ &\leq&
\delta^*_{\epsilon ^*}\mid w(\bar{x}) \mid + 2 \mid w(\bar{x}) \mid \epsilon
\end{eqnarray*}
as $ h\Big(\varphi\big(w(\bar{x})\big),id\Big)  \leq h\Big(\phi\big(w(\bar{x})\big),
\varphi\big(w(\bar{x})\big)\Big)
+  h\Big(\phi\big(w(\bar{x})\big), id\Big) $, then 
$$ 
h\Big(\phi\big(w(\bar{x})\big), \id\Big) \geq \alpha -\delta^*_{\epsilon ^*}\mid w(\bar{x}) \mid - 2 \mid w(\bar{x}) \mid 
\epsilon=\alpha_0
$$
\end{proof}
So, $\phi(p)\neq \id$ for sufficiently small $\epsilon>0$.
\end{proof}
\underline{\bf Acknowledgments}. The authors would like to thank E.
Gordon for useful discussions. The work was partially supported by
CONACyT grant SEP-25750, and PROMEP grant UASLP-CA-21.
\def\cprime{$'$} \def\cprime{$'$} \def\cprime{$'$}


IICO-UASLP, Av. Karakorum 1470, Lomas 4ta Session, San Luis Potosi,
SLP 7820 Mexico

e-mail:glebsky@cactus.iico.uaslp.mx

\end{document}